\documentclass[10pt]{article}

\usepackage{amsmath}
\usepackage{amsfonts}
\usepackage{amssymb}
\usepackage[english]{babel}
\usepackage{amsthm}

\newcommand{\mR}{\mathbb{R}}
\newcommand{\mC}{\mathbb{C}}

\newcommand{\mS}{\mathbb{S}}

\newcommand{\mcP}{\mathcal{P}}

\newcommand{\mcM}{\mathcal{M}}
\newcommand{\tmcM}{\widetilde{\mathcal{M}}}

\newcommand{\mcA}{\mathcal{A}}
\newcommand{\mcB}{\mathcal{B}}
\newcommand{\mcO}{\mathcal{O}}

\newcommand{\gf}{\mathfrak{f}}
\newcommand{\gfd}{\mathfrak{f}^{\dagger}}
\newcommand{\gsl}{\mathfrak{sl}}

\newcommand{\ol}{\overline}
\newcommand{\olz}{\ol{z}}

\newcommand{\uX}{\underline{X}}

\newcommand{\la}{\lambda}
\newcommand{\La}{\Lambda}

\newcommand{\uz}{\underline{z}}
\newcommand{\uzd}{\underline{z}^{\dagger}}
\newcommand{\tuz}{\widetilde{\uz}}
\newcommand{\tuzd}{\tuz^{\dagger}}

\newcommand{\upz}{\partial_{\uz}}
\newcommand{\tupz}{\widetilde{\upz}}
\newcommand{\upzd}{\partial_{\uz}^{\dagger}}
\newcommand{\tupzd}{\widetilde{\upzd}}

\newcommand{\tnz}{\vert\tuz\vert}
\newcommand{\p}{\partial}
\newcommand{\dirac}{\underline{\p}}

\newtheorem{theorem}{Theorem}

\newtheorem{remark}{Remark}
\newtheorem{corollary}{Corollary}

\begin{document}
\title{Gelfand-Tsetlin Bases of Orthogonal Polynomials in Hermitean Clifford Analysis}
\author{F.\ Brackx$^\ast$, H.\ De Schepper$^\ast$, R.\ L\'{a}vi\v{c}ka$^\ddagger$ \& V.\ Sou\v{c}ek$^\ddagger$}

\date{\small{$^\ast$ Clifford Research Group, Faculty of Engineering, Ghent University\\
Building S22, Galglaan 2, B-9000 Gent, Belgium\\
$^\ddagger$ Mathematical Institute, Faculty of Mathematics and Physics, Charles University\\
Sokolovsk\'a 83, 186 75 Praha, Czech Republic}}

\maketitle

\begin{abstract}
\noindent An explicit algorithmic construction is given for orthogonal bases for spaces of homogeneous polynomials, in the context of Hermitean Clifford analysis, which is a higher dimensional function theory centred around the simultaneous null solutions of two Hermitean conjugate complex Dirac operators. 
\end{abstract}



\section{Introduction}


This paper is devoted to the construction of orthogonal bases for spaces of Hermitean monogenic polynomials, and to establishing a recursive algorithm for this construction. The theory of Hermitean monogenic functions is one of the actual research topics in Clifford analysis, which, in its most basic form, is a higher dimensional generalization of holomorphic function theory in the complex plane, and a refinement of harmonic analysis, see e.g. \cite{bds, gilbert, dss, guerleb, ghs}. At the heart of this function theory lies the notion of a monogenic function, i.e.\ a Clifford algebra valued null solution of the Dirac operator $\dirac = \sum_{\alpha=1}^m e_{\alpha} \, \p_{X_{\alpha}}$, where $(e_1, \ldots, e_m)$ is an orthonormal basis of $\mR^m$ underlying the construction of the real Clifford algebra $\mR_{0,m}$. We refer to this setting as the Euclidean case, since the fundamental group leaving the Dirac operator $\dirac$ invariant is the orthogonal group $\mbox{O}(m;\mR)$, which is doubly covered by the Pin($m$) group of the Clifford algebra. In the books \cite{rocha,struppa} and the series of papers \cite{sabadini,eel,partI,partII,toulouse,eel2,howe} so--called Hermitean Clifford analysis recently emerged as a refinement of Euclidean Clifford analysis, where the considered functions now take their values in the complex Clifford algebra $\mC_m$ or in complex spinor space. Hermitean Clifford analysis is based on the introduction of an additional datum, a so--called complex structure $J$, inducing an associated Dirac operator $\dirac_J$; it then focusses on the simultaneous null solutions of both operators $\dirac$ and $\dirac_J$, called Hermitean monogenic functions. The corresponding function theory is still in full development, see also \cite{hehe, hermwav, mabo, mahi, cama, eel3}. It is worth mentioning that the traditional holomorphic functions of several complex variables are a special case of Hermitean monogenic functions.\\[-2mm]

The construction of orthogonal bases for spaces of null solutions of differential operators is much facilitated by the Gel'fand--Tsetlin approach, which will be used in this paper. The notion of Gel'fand--Tsetlin (GT) basis applies to finite dimensional irreducible modules over a classical Lie algebra (see \cite{GT,molev}); when this Lie algebra is realized explicitly, say as a subspace of null solutions of an invariant differential operator, then an algorithm for the construction of the GT-basis may be devised. Recently much effort has been put into the construction of orthogonal bases for spaces of monogenic polynomials in the framework of Euclidean Clifford analysis, to meet the needs for numerical calculations. Indeed, the basis polynomials, mostly called Fueter polynomials, appearing in the Taylor series expansion of monogenic functions, are not useful for that purpose since they are not orthogonal with respect to the natural $L_2$--inner product. Explicit constructions of orthogonal polynomial bases in the Euclidean Clifford analysis context were carried out in e.g. \cite{cacao,cagubo} in a direct analytic way starting from spherical harmonics, and in e.g. \cite{bock, lav,lasova,pvl} by the GT--approach.\\[-2mm]

The present paper describes in a systematic and detailed way the GT--construction of orthogonal bases for spaces of homogeneous Hermitean monogenic polynomials and in this way offers a full account of some results already announced in \cite{rhovla,rhoro,rhofb}. Moreover, special attention is paid  to the so--called Appell property of the constructed bases in complex dimension $2$, a property which is important for numerical applications. In Section 3 it is shown that, as is the case in Euclidean Clifford analysis, also in the present Hermitean case, two powerful tools, viz.\ the Cauchy-Kovalevskaya extension theorem and the Fischer decomposition theorem, are indispensable when implementing the Gel'fand--Tsetlin procedure. This GT--procedure together with the cited theorems enable devising the algorithm aimed at, which can be used in any dimension. Section 4 shows how the well-known orthogonal bases of polynomials in the complex plane fit into the general scheme. Sections 5 and 6 treat the exceptional cases where Hermitean monogenicity reduces to (anti--)holomorphy in several complex variables. In Section 7, the algorithm is applied to full extent in order to explicitly construct the orthogonal bases in complex dimension $n=2$. In this particular dimension, the constructed bases show a special and important property. Each of the four complex derivatives maps any basis polynomial of a certain given bidegree of homogeneity, to a simple multiple of another basis element of the corresponding lower bi-homogeneity degree. It means that all four derivatives are represented, in this particular basis, by a matrix with a simple structure; details are contained in Section 8. Finally, in the last section, using the fact that any solution of the Hermitean monogenicity system of equations is automatically monogenic, we are able to construct a suitable basis for spaces of monogenic polynomials of a given degree of homogeneity  in terms of the orthogonal bases of Hermitean monogenic polynomials constructed in this paper. To make the paper self--contained an introductory section on Clifford analysis is included.


\section{Preliminaries on Clifford analysis}


For a detailed description of the structure of Clifford algebras we refer to e.g.\ \cite{port}. Here we only recall the necessary basic notions. The real Clifford algebra $\mathbb{R}_{0,m}$ is constructed over the vector space $\mathbb{R}^{0,m}$ endowed with a non--degenerate quadratic form of signature $(0,m)$ and generated by the orthonormal basis $(e_1,\ldots,e_m)$. The non--commutative Clifford or geometric multiplication in $\mathbb{R}_{0,m}$ is governed by the rules 
\begin{equation}\label{multirules}
e_{\alpha} e_{\beta} + e_{\beta} e_{\alpha} = -2 \delta_{\alpha \beta} \ \ , \ \ \alpha,\beta = 1,\ldots ,m
\end{equation}
As a basis for $\mathbb{R}_{0,m}$ one takes for any set $A=\lbrace j_1,\ldots,j_h \rbrace \subset \lbrace 1,\ldots,m \rbrace$ the element $e_A = e_{j_1} \ldots e_{j_h}$, with $1\leq j_1<j_2<\cdots < j_h \leq m$, together with $e_{\emptyset}=1$, the identity element. Any Clifford number $a$ in $\mathbb{R}_{0,m}$ may thus be written as $a = \sum_{A} e_A a_A$, $a_A \in \mathbb{R}$, or still as $a = \sum_{k=0}^m \lbrack a \rbrack_k$, where $\lbrack a \rbrack_k = \sum_{|A|=k} e_A a_A$ is the so--called $k$--vector part of $a$. Euclidean space $\mathbb{R}^{0,m}$ is embedded in $\mathbb{R}_{0,m}$ by identifying $(X_1,\ldots,X_m)$ with the Clifford vector $\uX = \sum_{\alpha=1}^m e_{\alpha}\, X_{\alpha}$, for which it holds that $\uX^2 = - |\uX|^2$. The vector valued first order differential operator $\dirac = \sum_{\alpha=1}^m e_{\alpha}\, \p_{X_{\alpha}}$, called Dirac operator, is the Fourier or Fischer dual of the Clifford variable $\uX$. It is this operator which underlies the notion of monogenicity of a function, a notion which is the higher dimensional counterpart of holomorphy in the complex plane. More explicitly, a function $f$ defined and continuously differentiable in an open region $\Omega$ of $\mathbb{R}^{m}$ and taking values in (a subspace of) the Clifford algebra $\mathbb{R}_{0,m}$, is called (left) monogenic in $\Omega$ if $\dirac\lbrack f \rbrack = 0$ in $\Omega$. As the Dirac operator factorizes the Laplacian: $\Delta_m = - \dirac^2$, monogenicity can be regarded as a refinement of harmonicity. The Dirac operator being rotationally invariant, this framework is usually referred to as Euclidean Clifford analysis.\\[-2mm]

When allowing for complex constants, the generators $(e_1,\ldots, e_{m})$, still satisfying (\ref{multirules}), produce the complex Clifford algebra $\mathbb{C}_{m} = \mathbb{R}_{0,m} \oplus i\, \mathbb{R}_{0,m}$. Any complex Clifford number $\lambda \in \mathbb{C}_{m}$ may thus be written as $\lambda = a + i b$, $a,b \in \mathbb{R}_{0,m}$, leading to the definition of the Hermitean conjugation $\lambda^{\dagger} = (a +i b)^{\dagger} = \overline{a} - i \overline{b}$, where the bar notation stands for the Clifford conjugation in $\mathbb{R}_{0,m}$, i.e. the main anti--involution for which $\overline{e}_{\alpha} = -e_{\alpha}$, $\alpha=1, \ldots,m$. This Hermitean conjugation leads to a Hermitean inner product on $\mathbb{C}_{m}$ given by $(\lambda,\mu) = \lbrack \lambda^{\dagger} \mu \rbrack_0$ and its associated norm $|\lambda| = \sqrt{ \lbrack \lambda^{\dagger} \lambda \rbrack_0} = ( \sum_A |\lambda_A|^2 )^{1/2}$. This is the framework for Hermitean Clifford analysis, which emerges from Euclidean Clifford analysis by introducing an additional datum, a so--called complex structure, i.e.\ an $\mbox{SO}(m;\mR)$--element $J$ with $J^2=-\mathbf{1}$ (see \cite{partI,partII}), forcing the dimension to be even; from now on we put $m=2n$. Usually $J$ is chosen to act upon the generators of $\mC_{2n}$ as $J[e_j] = -e_{n+j}$ and $J[e_{n+j}] = e_j$, $j=1,\ldots,n$. By means of the projection operators $\pm \frac{1}{2}(\mathbf{1} \pm iJ)$ associated to $J$, first the Witt basis elements $(\gf_j,\gf_j^{\dagger})^n_{j=1}$ for $\mC_{2n}$ are obtained: 
\begin{eqnarray*}
\gf_j & = & \phantom{-}\frac{1}{2} (\mathbf{1} + iJ) [e_j] = \phantom{-}\frac{1}{2} (e_j - i \, e_{n+j}), \qquad j=1,\ldots,n \\
\gf_j^\dagger & = & -\frac{1}{2} (\mathbf{1} - iJ)[e_j] =  -\frac{1}{2} (e_j + i\, e_{n+j}), \qquad j=1,\ldots,n
\end{eqnarray*}
The Witt basis elements satisfy the respective Grassmann and duality identities
$$
\gf_j \gf_k + \gf_k \gf_j = \gf_j^{\dagger} \gf_k^{\dagger} + \gf_k^{\dagger} \gf_j^{\dagger} = 0, \quad
\gf_j \gf_k^{\dagger} + \gf_k^{\dagger} \gf_j = \delta_{jk}, \ \ j,k=1,\ldots, n
$$
whence they are isotropic: $(\gf_j)^2 = 0, (\gfd_j)^2 = 0, j=0,\ldots,n$.
Next, a vector in $\mR^{0,2n}$ is now denoted by $(x_1, \ldots , x_n, y_1 , \ldots , y_n)$ and identified with the Clifford vector $\uX  =  \sum_{j=1}^n (e_j\, x_j + e_{n+j}\, y_j)$, producing, by projection, the Hermitean Clifford variables $\uz$ and $\uzd$:
$$
\uz = \frac{1}{2} (\mathbf{1} + iJ) [\uX] = \sum_{j=1}^n \gf_j\, z_j, \quad
\uzd = -\frac{1}{2} (\mathbf{1} - iJ)[\uX] = \sum_{j=1}^n \gf_j^{\dagger}\, z_j^c
$$
where complex variables $z_j = x_j + i y_j$ have been introduced, with complex conjugates $z_j^c = x_j - i y_j$, $j=1,\ldots,n$. Finally, the Euclidean Dirac operator $\dirac$ gives rise, in the same manner, to the Hermitean Dirac operators $\upz$ and $\upzd$:
$$
\upzd = \frac{1}{4} (\mathbf{1} + iJ) [\dirac] = \sum_{j=1}^n \gf_j\, \p_{z^c_j}, \quad
\upz = -\frac{1}{4} (\mathbf{1} - iJ)[\dirac] = \sum_{j=1}^n \gf_j^{\dagger}\, \p_{z_j}  
$$
involving the Cauchy--Riemann operators $\p_{\olz_j} = \frac{1}{2} (\p_{x_j} + i \p_{y_j})$  and their complex conjugates $\p_{z_j} = \frac{1}{2} (\p_{x_j} - i \p_{y_j})$  in the $z_j$--planes, $j=1,\ldots,n$. Observe that Hermitean vector variables and Dirac operators are isotropic, i.e. $\uz^2 = (\uzd)^2 = 0$ and $(\upz)^2 = (\upzd)^2 = 0$, whence the Laplacian allows for the decomposition and factorization
$$
\Delta_{2n} = 4(\upz \upzd + \upzd \upz) = 4(\upz + \upzd)^2 = -4(\upzd - \upz)^2
$$
while dually
\begin{displaymath}
-(\uz-\uzd)^2 = (\uz + \uzd)^2 = \uz\, \uzd + \uzd \uz = |\uz|^2 = |\uzd|^2 = |\uX|^2
\end{displaymath}

\vspace*{3mm}
We consider functions with values in an irreducible representation $\mS_n$ of $\mC_{2n}$, called spinor space, which is realized within $\mC_{2n}$ using a primitive idempotent $I = I_1 \ldots I_n$, with $I_j = \gf_j \gf_j^{\dagger}$, $j=1,\ldots,n$. With that choice, $\gf_jI=0, j=1,\ldots,n$, and so $\mS_n \equiv \mC_{2n} I \cong {\bigwedge}_n^{\dagger} I$, where ${\bigwedge}_n^{\dagger}$ denotes the Grassmann algebra generated by the $\gf_j^\dagger$'s. Hence $\mS_n$ decomposes into homogeneous parts as 
$$
\mS_n = \bigoplus\limits_{r=0}^n \ \mS_n^{(r)} = \bigoplus\limits_{r=0}^n \ ( {\bigwedge}_n^{\dagger} )^{(r)} I
$$
with $( {\bigwedge}_n^\dagger )^{(r)} = \mbox{span}_{\mC} ( \gfd_{k_1} \wedge \gfd_{k_2} \wedge \cdots \wedge \gfd_{k_r} : \{ k_1,\ldots,k_r \} \subset \{1,\ldots,n  \} )$.\\[-2mm]

A continuously differentiable function $g$ in an open region $\Omega$ of $\mathbb{R}^{2n}$ with values in (a subspace of) the complex Clifford algebra $\mathbb{C}_{2n}$ then is called (left) Hermitean monogenic (or h--monogenic) in $\Omega$ if and only if it satisfies in $\Omega$ the system $\upz\, g = 0 = \upzd\, g$, or, equivalently, the system $\dirac \, g = \dirac_J \, g$, with $\dirac_J = J[\dirac]$. A major difference between Hermitean and Euclidean Clifford analysis concerns the underlying group invariance. Where the Euclidean Dirac operator $\dirac$ is invariant under the action of SO$(m)$, the system invariance of $(\upz,\upzd)$ breaks down to the group U$(n)$, see e.g.\ \cite{partI,partII}. For this reason U$(n)$ will play a fundamental role in our construction of orthogonal bases of spaces of Hermitean monogenic polynomials.\\[-2mm]

The spaces of homogeneous polynomials on $\mC^{n}$ with bidegree of homogeneity $(a,b)$ in $(\uz,\uzd)$ and taking values in $\mS_n^{(r)}$, will be denoted by $\mcP_{a,b}^{(r)}(\mC^n)$. By $\mcM_{a,b}(\mC^n)$ we denote the space of Hermitean monogenic polynomials of bidegree $(a,b)$ in $(\uz,\uzd)$, and by $\mcM_{a,b}^{(r)}(\mC^n)$ its subspace with values in $\mS_n^{(r)}$.


\section{The Gel'fand--Tsetlin procedure for U($n$)-modules}


As was already pointed out in the foregoing section, the fundamental group of Hermitean Clifford analysis in $\mC^n$ is the unitary group U$(n)$. The construction of a  Gel'fand--Tsetlin (GT) basis for an irreducible U$(n)$--module first requires a chain of subgroups U$(n)\supset$ U$(n-1)\supset\ldots\supset \; $U$(1)$. We choose these embeddings in such a way that in each step the last variable is preserved by the corresponding subgroup. Next we need the branching rules for  U$(n)$--modules, which are expressed in terms of  the highest weights of the corresponding irreducible modules. Irreducible U$(n)$--modules are classified according to their highest weight $\la=(\la_1,\ldots,\la_n)$, where the integers $\la_i$ satisfy the traditional condition $\la_1\geq\la_2\geq\ldots\geq\la_n$. The branching rules then are the following, see e.g.\ \cite{good}.

\begin{theorem} 
\label{BR}
When restricting to U$(n-1)$, an irreducible U$(n)$-module $V_\la$ with highest weight $\la$ decomposes as
$$
V_\la=\bigoplus_{\mu \prec \la} \ V_\mu
$$
where each of the summands $V_\mu$ is an irreducible U$(n-1)$--module with highest weight $\mu = (\mu_1, \mu_2, \ldots, \mu_{n-1})$; this multiplicity free direct sum is taken over all possible highest weights $\mu$ such that $\la \succ \mu$, i.e. $\la_1\geq\mu_1\geq\la_2\geq\ldots\geq\la_{n-1}\geq\mu_{n-1}\geq\la_n$.	If, moreover, $V_{\la}$ is endowed with a U$(n)$--invariant inner product, then this decomposition is orthogonal w.r.t. that inner product.
\end{theorem}

\noindent Proceeding by induction Theorem 1 generates the GT--basis for $V_\la$; this result reads as follows.

\begin{theorem}
With respect to the chain of subgroups U$(n)\supset$ U$(n-1)\supset\ldots\supset \; $U$(1)$, the irreducible U$(n)$--module $V_{\la^{(n)}}$ with highest weight $\la^{(n)} = (\la^{(n)}_1, \ldots, \la^{(n)}_n)$ decomposes into a direct sum of one-dimensional subspaces 
$$
V_{\la^{(n)}}  =  \bigoplus_{\La} \ V_\La
$$
taken over all possible GT--patterns $\La = \left(\la^{(n)},\ldots,\la^{(1)} \right)$ such that $\la^{(n)} \succ \la^{(n-1)} \succ \ldots \succ \la^{(1)}$, with $\la^{(j)}=(\la^{(j)}_1,\ldots,\la^{(j)}_j)$.
If moreover $V_{\la^{(n)}}$ is endowed with a U$(n)$--invariant inner product, then the above decomposition is orthogonal w.r.t. that inner product.
\end{theorem}

\noindent Note that each module $V_\La$ is uniquely determined by its GT--pattern $\La = \left(\la^{(n)},\ldots,\la^{(1)} \right)$, which in its turn depends uniquely on the chosen chain of subgroups of U$(n)$.\\[-2mm]

This Gel'fand--Tsetlin procedure will now be used to establish an induction algorithm allowing for the explicit construction of orthogonal bases for spaces of spherical Hermitean monogenics, defined in $\mC^n$ and taking values in the homogeneous parts $\mS_n^{(r)}$ of spinor space $\mS_n$. To that end it is quite necessary to explicitly describe the branching rule of Theorem 1. Let us have a closer look at the induction step from $(n-1)$ to $n$. By the induction hypothesis it is assumed that the GT--bases are known, in dimension $(n-1)$, for all spaces $\tmcM_{a,b}^{(r)}(\mC^{n-1})$ of Hermitean monogenic polynomials in the variables $(\tuz,\tuzd) = (z_1, \ldots, z_{n-1},\olz_1, \ldots, \olz_{n-1} )$, which are homogeneous of bidegree $(a,b)$ and take their values in $\mS_{n-1}^{(r)}, r=0, \ldots, n-1$. Now consider the space $\mcM_{a,b}^{(r)}(\mC^n)$ with $a$, $b$ and $r$ fixed. We will consider the general case $0<r<n$, since for $r=0$ and $r=n$ there are small variations in the form of the highest weights which do not fit into the general approach. As a matter of fact the exceptional cases $r=0$ and $r=n$ are quite easily treated since in those cases the notion of Hermitean monogenicity is nothing else but (anti-)holomorphy in several complex variables; those two particular cases will be treated in Sections 5 and 6. The highest weight of $\mcM_{a,b}^{(r)}(\mC^n)$ is $\la^{(n)} = (a+1, 1, \ldots, 1, 0, \ldots, 0, -b)$ with the last component $1$ at the $r$-th place. We denote this weight shortly by $\la^{(n)} =[a,-b]_r$. From Theorem 1 we know that the irreducible U$(n)$-module $\mcM_{a,b}^{(r)}(\mC^n)$ may be decomposed as
$$
\mcM_{a,b}^{(r)}(\mC^n) = \bigoplus_{\mu \prec \la^{(n)}} \ V_{\mu}
$$
where the $V_{\mu}$ are irreducible U$(n-1)$-modules with highest weights $\mu \prec \la^{(n)} = [a,-b]_r$. Those highest weights are clearly of the form $[j,k]_r$ and $[j,k]_{r-1}$ with $j=0,\ldots,a$ and $k=-b,\ldots,0$. These may be interpreted as the co-ordinates of the lattice points in two rectangular grids, one for the spinor homogeneity degree $r$ and one for $(r-1)$, with vertices $(0,0)$, $(a,0)$, $(0,-b)$ and $(a,-b)$ (see Figures 1 and 2). This means that in the above decomposition $2(a+1)(b+1)$ U$(n-1)$-modules of polynomials of different homogeneity bidegrees $(j,-k)$, $j=0,\ldots,a$, $k=-b,\ldots,0$ are needed. Moreover these U$(n-1)$-modules have to be embedded into $\mcM_{a,b}^{(r)}(\mC^n)$ while preserving the U$(n-1)$-invariance. We will now show that Hermitean Clifford analysis provides the adequate tools to achieve this goal.\\[-2mm]

Embedding homogeneous polynomials into a space of homogeneous Hermitean monogenic polynomials such as $\mcM_{a,b}^{(r)}(\mC^n)$ while keeping U$(n-1)$-invariance, is exactly the context of the so-called Cauchy-Kovalevskaya (CK) extension principle, which, originally, is about the existence and uniqueness of solutions of initial value problems for PDE (see \cite{cooke} for a historic overview). The classical idea of the CK extension is to characterize solutions of suitable (systems of) PDE's by their restriction, sometimes together with the restrictions of a number of their derivatives, to a submanifold of real codimension $1$. A well--known typical result in this context concerns a real-analytic function $f(x)$ in a neighbourhood of the origin on the real axis, for which there exists the CK-operator $\exp{\left(iy\frac{d}{dx}\right)}$, such that 
$$
F(x) = \exp{\left(iy\frac{d}{dx}\right)} \left[ f(x) \right] = \sum_{k=0}^{\infty} \ \frac{1}{k!} (iy)^k f^{(k)}(x)
$$
is holomorphic in a neigbourhood of the origin in the complex plane, and its restriction to the real axis is exactly $f(x)$. Even so, in Euclidean Clifford analysis, for a real-analytic function $f(\uX)$ in a neighbourhood of the origin in $\mR^m$, there is the CK-operator $\exp{\left(X_m e_m \widetilde{\dirac}\right)}$ such that
$$
F(\uX) = \exp{\left(X_m e_m \widetilde{\dirac}\right)} \left[ f(\uX) \right] = \sum_{k=0}^{\infty} \ \frac{1}{k!} \ X_m^k (e_m \widetilde{\dirac})^k \left[ f(\uX) \right]
$$
is monogenic in a neighbourhood of the origin in $\mR^{m+1}$ (see \cite{bds}) and its restriction to $\mR^m$ is $f(\uX)$.\\[-2mm]

Recently, in \cite{CK}, we have obtained a CK-extension theorem for Hermitean monogenic functions by restricting the null solutions of the Hermitean Dirac operators, and a number of their derivatives, to the vector subspace $\mC^{n-1}$ of complex codimension $1$. To that end we single out the variables $(z_n,\olz_n)$ and consider restrictions to $\{z_n = 0 = \olz_n \}$ identified with $\mC^{n-1}$. The value space $\mS_n^{(r)} = \left(\bigwedge_n^{\dagger}\right)^{(r)}\left(\gfd_1, \ldots, \gfd_n\right) I$ is then split as
$$
\mS_n^{(r)} = \left({\bigwedge}_{n-1}^{\dagger}\right)^{(r)}\left(\gfd_1, \ldots, \gfd_{n-1}\right) I \oplus \gfd_n \ \left({\bigwedge}_{n-1}^{\dagger}\right)^{(r-1)}\left(\gfd_1, \ldots, \gfd_{n-1}\right) I
$$
and functions are split accordingly as $F = F^0 + \gfd_n F^1$, where $F^0$ takes its values in $(\bigwedge_{n-1}^{\dagger})^{(r)} I$, while $F^1$ takes its values in $(\bigwedge_{n-1}^{\dagger})^{(r-1)} I$. In the same order of ideas the Hermitean variables and the Hermitean Dirac operators are split as
$$
\uz = \tuz + \gf_n z_n \quad , \quad \uzd = \tuzd + \gfd_n \olz_n
$$
and
$$
\upz = \tupz + \gfd_n \p_{z_n} \quad , \quad \upzd = \tupzd + \gf_n \p_{\olz_n}
$$

The CK-extension theorem of Hermitean Clifford analysis then reads as follows (see \cite{CK}).

\begin{theorem}
Given homogeneous polynomials $p_{a,b-j}^0$, $j=0,\ldots,b$ and $p_{a-i,b}^1$, $i=0,\ldots,a$, respectively taking values in $({\bigwedge}_{n-1}^{\dagger})^{(r)} (\gfd_1, \ldots, \gfd_{n-1}) \, I$ and in $({\bigwedge}_{n-1}^{\dagger})^{(r-1)} (\gfd_1, \ldots, \gfd_{n-1}) \, I$, and satisfying the respective compatibility conditions
\begin{eqnarray*}
&& \tupz p_{a,b}^0 = 0, \ \tupz p_{a,b-1}^0 = 0, \ \ldots, \ \tupz p_{a,0}^0 = 0 \quad (r < n-1) \\
&& \tupzd p_{a,b}^1 = 0, \ \tupzd p_{a-1,b}^1 = 0, \  \ldots, \ \tupzd p_{0,b}^1 = 0 \quad (r > 1)
\end{eqnarray*}
there exists a unique Hermitean monogenic homogeneous polynomial $M_{a,b}$ such that
\begin{itemize}
\item[(i)] $\displaystyle{M_{a,b}|_{\mC^{n-1}} = p_{a,b} = p_{a,b}^0 + \gfd_n p_{a,b}^1}$;
\item[(ii)] $\displaystyle{\frac{\p^j}{\p_{z_n^c}^j} M_{a,b}|_{\mC^{n-1}} = p_{a,b-j} = p_{a,b-j}^0 - \gfd_n \tupzd p_{a,b-j+1}^0}$, $j=1,\ldots,b$;
\item[(iii)] $\displaystyle{\frac{\p^i}{\p_{z_n}^i} M_{a,b}|_{\mC^{n-1}} = p_{a-i,b} =  \tupz p_{a-i+1,b}^1 + \gfd_n p_{a-i,b}^1}$,  $i=1,\ldots,a$.
\end{itemize}
This CK--extension $M_{a,b}$ is given by
$$
M_{a,b} = \sum_{j=0}^b \, M_{a,b-j}^0 + \sum_{i=0}^a \, M_{a-i,b}^1
$$
where
$$
M_{a,b-j}^0 = {(\olz_n)}^j \ \sum_{k=0}^{\min{(2a+1,2b-2j)}} \frac{1}{\lfloor \frac{k}{2} \rfloor!}\frac{1}{\lfloor \frac{k+1}{2} + j \rfloor!} \left( z_n \, \tupz \, \gf_n \, +  \, \olz_n \, \tupzd \, \gfd_n \right)^k \left[ \, p_{a,b-j}^0 \, \right]
$$
and
$$
M_{a-i,b}^1 = z_n^i \ \sum_{k=0}^{\min{(2a-2i,2b+1)}} \frac{1}{\lfloor \frac{k}{2} \rfloor!}\frac{1}{\lfloor \frac{k+1}{2} + i \rfloor!} \left( z_n \, \tupz \, \gf_n +  \olz_n \, \tupzd \gfd_n \, \right)^k \left[\, \gfd_n \, p_{a-i,b}^1 \, \right]
$$
\end{theorem}

\noindent In fact this means that there are $a+b+2$ spaces of initial data, viz.
$$
\mcA_{a,b-j}^{(r)} = \{ p_{a,b-j}^0 \in \mcP_{a,b-j}^{(r)}(\mC^{n-1}) : \tupz \, p_{a,b-j}^0 = 0 \} \quad , \quad j=0,\ldots,b
$$
and
$$
\mcB_{a-i,b}^{(r)} = \{ \gfd_n \, p_{a-i,b}^1 :   p_{a-i,b}^1 \in \mcP_{a-i,b}^{(r-1)}(\mC^{n-1}) \quad {\rm and} \quad \tupzd \, p_{a-i,b}^1 = 0 \} \quad , \quad i=0,\ldots,a
$$
which by the CK-extension procedure are mapped isomorphically to subspaces of $\mcM_{a,b}^{(r)}(\mC^n)$. In this way the whole of $\mcM_{a,b}^{(r)}(\mC^n)$ is recovered, which is confirmed by counting dimensions. Indeed, we know from \cite{fischer,CK} that
\begin{eqnarray*}
x_{a,b,r} & = & {\rm dim}(\mcA_{a,b-j}^{(r)}) \ = \ \frac{r}{a+r} {n-1 \choose r} {a+n-1 \choose a} {b+n-2 \choose b}  \\
y_{a,b,r} & = & {\rm dim}(\mcB_{a-i,b}^{(r)}) \ = \ \frac{r}{b+n-r} {n-1 \choose r} {a+n-2 \choose a} {b+n-1 \choose b}
\end{eqnarray*}
while
$$
m_{a,b}^{(r)} = {\rm dim}(\mcM_{a,b}^{(r)}(\mC^n)) =   \frac{r(a+b+n)}{(a+r)(b+n-r)} {n-1 \choose r} {a+n-1 \choose a} {b+n-1 \choose b}
$$
and it may be readily verified that
$$
\sum_{j=0}^b \ x_{a,j,r} + \sum_{i=0}^a \ y_{i,b,r} = m_{a,b}^{(r)} 
$$
It is important to note that the dimension $m_{a,b}^{(r)} = {\rm dim}(\mcM_{a,b}^{(r)}(\mC^n))$ was determined in \cite{fischer} independently of the CK-extension context, even by two approaches, one of them involving the decomposition of harmonic homogeneous polynomials in terms of Hermitean monogenic ones, the other being based on the Weyl dimension formula (see \cite{good}, p.382).\\[-2mm]

As the CK-extension operator commutes with the action of U$(n-1)$, we will obtain, by the CK-extension, a decomposition of $\mcM_{a,b}^{(r)}(\mC^n)$ into a direct sum of U$(n-1)$-invariant subspaces, if at least we are able to decompose the spaces of initial polynomials $\mcA_{a,b-j}^{(r)}$ and $\mcB_{a-i,b}^{(r)}$ into U$(n-1)$-irreducibles. This kind of decomposition is the content of the so-called Fischer decomposition theorem, which is, quite naturally, very well-known for spherical harmonics, but also exists in the Euclidean and the Hermitean Clifford analysis setting (see \cite{dss,damdeef,howe,fischer}). In particular the kernels of the Hermitean Dirac operators may be decomposed in terms of spherical Hermitean monogenics, a result which reads, in a form adapted to the present situation, as follows (see \cite{rhoro}).

\begin{theorem}
(i) Under the action of U$(n-1)$ the space $\mcA_{a,b}^{(r)} = \{ p_{a,b}^0 \in \mcP_{a,b}^{(r)}(\mC^{n-1}) : \tupz p_{a,b}^0 = 0 \}$ has the multiplicity free irreducible decomposition
\begin{eqnarray*}
\mcA_{a,b}^{(r)} &=& \tmcM_{a,b}^{(r)} \oplus \bigoplus_{k=0}^{\min(a,b-1)} \tnz^{2k} \tuzd \tmcM_{a-k,b-k-1}^{(r-1)}\\
& & \phantom{\tmcM_{a,b}^{(r)}} \oplus \bigoplus_{k=0}^{\min(a-1,b-1)} \tnz^{2k} \left( \tuzd \tuz + \frac{a-k-1+r}{a+r} \, \tuz \tuzd \right) \tmcM_{a-k-1,b-k-1}^{(r)}
\end{eqnarray*}
(ii) Under the action of U$(n-1)$ the space $\mcB_{a,b}^{(r)} = \{ \gfd_n p_{a,b}^1 :   p_{a,b}^1 \in \mcP_{a,b}^{(r-1)}(\mC^{n-1}) \; {\rm and} \; \tupzd p_{a,b}^1 = 0 \}$ has the multiplicity free irreducible decomposition
\begin{eqnarray*}
\mcB_{a,b}^{(r)} &=& \gfd_n (  \tmcM_{a,b}^{(r-1)} \oplus \bigoplus_{k=0}^{\min(a-1,b)} \tnz^{2k} \tuz \tmcM_{a-k-1,b-k}^{(r)}\\
& &  \phantom{\gfd_n (  \tmcM_{a,b}^{(r-1)}}\oplus \bigoplus_{k=0}^{\min(a-1,b-1)} \tnz^{2k} \left( \tuz \tuzd + \frac{b-k-1+n-r}{b+n-r} \, \tuzd \tuz \right) \tmcM_{a-k-1,b-k-1}^{(r-1)})
\end{eqnarray*}
\end{theorem}

Combining the CK-extension and Fischer decomposition theorems now leads to an inductive algorithm for the Gel'fand-Tsetlin construction of a basis for $\mcM_{a,b}^{(r)}(\mC^{n})$, which will be orthogonal w.r.t. the natural inner product:\\[-2mm]

\noindent (i) {\em step $0$ (induction hypothesis)} \\
\noindent keeping $a,b,r,n$ fixed, it is assumed that the orthogonal bases of each of the spaces $\mcM_{k,l}^{(s)}(\mC^{n-1})$, $0 \leq s \leq n-1$, $0 \leq k \leq a, 0 \leq l \leq b$ are known; \\[-2mm]

\noindent (ii) {\em step $1$} \\
\noindent by means of Theorem 4, the initial data spaces $\mcA_{a,b-j}^{(r)}$, $j=0,\ldots,b$ and $\mcB_{a-i,b}^{(r)}$, $i=0,\ldots,a$ are decomposed in terms of U$(n-1)$-irreducibles, which are shifted versions of the $\mcM_{k,l}^{(s)}(\mC^{n-1})$, obtaining, by means of the induction hypothesis, their basis polynomials;\\[-2mm]

\noindent (iii) {\em step $2$} \\
\noindent by means of Theorem 3, construct the CK-extension of each of the initial basis polynomials obtained in step $1$; these CK--extensions taken together form an orthogonal basis of $\mcM_{a,b}^{(r)}(\mC^n)$.\\[-2mm]

Now we will show that this algorithmic construction matches the decomposition of $\mcM_{a,b}^{(r)}(\mC^n)$ into U$(n-1)$-irreducibles $V_{\mu}$ with $\mu \prec \la^{(n)} = (a+1, 1, \ldots, 1, 0, \ldots, 0,-b) = [a,-b]_r$, as predicted by Theorem 1. In the Fischer decomposition of the initial data space $\mcA_{a,b}^{(r)}$ into U$(n-1)$-irreducibles, the set of highest weights of the components with values in $\mS_{n-1}^{(r)}$, is the string $[a,-b]_r, [a-1,-b+1]_r, [a-2,-b+2]_r, \ldots$ ending either with $[a-b,0]_r$ or with $[0,-b+a]_r$. The components with values in $\mS_{n-1}^{(r-1)}$ give rise to the string of highest weights $[a,-b+1]_{r-1}, [a-1,-b+2]_{r-1}, \ldots$ ending either with $[a-b+1,0]_{r-1}$ or with $[0,-b+a+1]_{r-1}$. Those two strings correspond to an anti-diagonal segment in the lattices introduced above (see Figure 1 for the case $a>b$ and Figure 2 for the case $a<b$). Similarly, for the initial data space $\mcA_{a,b-1}^{(r)}$, we obtain the strings of highest weights $[a,-b+1]_r, [a-1,-b+2]_r, \ldots$, ending either with $[a-b+1,0]_r$ or $[0,-b+a+1]_r$, and $[a,-b+2]_{r-1}, [a-1,-b+3]_{r-1}, \ldots$, ending either with $[a-b+2,0]_{r-1}$ or with $[0,-b+a+2]_{r-1}$, which, again, correspond to an anti-diagonal segment in each of the two lattices. Continuing in the same way, we obtain the strings $[a,-1]_r, [a-1,0]_r$ and $[a,0]_{r-1}$ for $\mcA_{a,1}^{(r)}$ and finally the string $[a,0]_r$ for $\mcA_{a,0}^{(r)}$. Starting with the initial data space $\mcB_{a,b}^{(r)}$, we obtain the strings of highest weights $[a,-b]_{r-1}, [a-1,-b+1]_{r-1}, \ldots$, ending either with $[a-b,0]_{r-1}$ or $[0,-b+a]_{r-1}$, and $[a-1,-b]_{r}, [a-2,-b+1]_r, \ldots$, ending either with $[a-b-1,0]_r$ or $[0,-b+a-1]_r$. The last string $[0,-b]_{r-1}$ is obtained from $\mcB_{0,b}^{(r)}$, and finally it is found that all predicted highest weights $\mu$ represented by all the lattice points in the two rectangular grids, are recovered by our construction. In fact, in this way we have given, by the sole use of Hermitean Clifford analysis tools, an independent proof of the branching rules for the representations which are realized as spaces of Hermitean monogenic polynomials.\\[-2mm]

Finally note that on each irreducible representation of U$(n)$ there always exists an invariant inner product, uniquely determined up to a constant. On $\mcM_{a,b}^{(r)}$ two well--known realizations of this inner product  are the $L_2$ inner product and the Fischer inner product, this last one being given for $P(\uz,\uzd)$ and $Q(\uz,\uzd)$ in $\mcM_{a,b}^{(r)}$ by 
$$
\langle P , Q \rangle = P^{\dagger}(\upzd,\upz)\left[Q(\uz,\uzd)\right]
$$

\unitlength 0.65cm
\begin{figure}[h]
\begin{picture}(22,8.5)(-0.5,0)
\put(-0.15,1.85){$\bullet$}
\put(0,2.0){\line(1,0){10}}
\put(5.85,1.85){$\bullet$}
\put(8.85,1.85){$\bullet$}
\put(9.85,1.85){$\bullet$}
\put(10,2.0){\line(0,1){6}}
\put(9.85,2.85){$\bullet$}
\put(9.85,7.85){$\bullet$}
\put(10,8.0){\line(-1,0){10}}
\put(4.85,7.85){$\bullet$}
\put(3.85,7.85){$\bullet$}
\put(2.85,7.85){$\bullet$}
\put(-0.15,7.85){$\bullet$}
\put(0,8.0){\line(0,-1){6}}
\put(0,8.0){\line(1,-1){6}}
\put(3,8.0){\line(1,-1){6}}
\put(4,8.0){\line(1,-1){6}}
\put(5,8.0){\line(1,-1){5}}
\put(-0.5,8.3){\footnotesize{$(0,0)$}}
\put(-0.5,1.5){\footnotesize{$(0,$-$b)$}}
\put(5.5,1.5){\footnotesize{$(b,$-$b)$}}
\put(9.5,1.5){\footnotesize{$(a,$-$b)$}}
\put(7.9,1.5){\footnotesize{$(a$-$1,$-$b)$}}
\put(9.5,8.3){\footnotesize{$(a,0)$}}
\put(3.2,8.3){\footnotesize{$(a$-$b,0)$}}
\put(4.6,8.3){\footnotesize{$(a$-$b$+$1,0)$}}
\put(1.4,8.3){\footnotesize{$(a$-$b$-$1,0)$}}
\put(9.3,3.3){\footnotesize{$(a,$-$b$+$1)$}}
\put(4.0,0.0){$\mS_{n-1}^{(r)}$ valued}
\put(16,0.0){$\mS_{n-1}^{(r-1)}$ valued}
\put(11.85,1.85){$\bullet$}
\put(12,2.0){\line(1,0){10}}
\put(12.85,1.85){$\bullet$}
\put(21.85,1.85){$\bullet$}
\put(22,2.0){\line(0,1){6}}
\put(21.85,1.85){$\bullet$}
\put(21.85,3.85){$\bullet$}
\put(21.85,7.85){$\bullet$}
\put(22,8.0){\line(-1,0){10}}
\put(17.85,7.85){$\bullet$}
\put(16.85,7.85){$\bullet$}
\put(15.85,7.85){$\bullet$}
\put(11.85,7.85){$\bullet$}
\put(12,8.0){\line(0,-1){6}}
\put(11.85,7.85){$\bullet$}
\put(12,3.0){\line(1,-1){1}}
\put(16,8.0){\line(1,-1){6}}
\put(17,8.0){\line(1,-1){5}}
\put(18,8.0){\line(1,-1){4}}
\put(11.85,2.85){$\bullet$}
\put(21.85,2.85){$\bullet$}
\put(11.5,8.3){\footnotesize{$(0,0)$}}
\put(11.5,1.5){\footnotesize{$(0,$-$b)$}}
\put(12.7,1.5){\footnotesize{$(1,$-$b)$}}
\put(21.5,8.3){\footnotesize{$(a,0)$}}
\put(21.5,1.5){\footnotesize{$(a,$-$b)$}}
\put(11.3,3.3){\footnotesize{$(0,$-$b$+$1)$}}
\put(21.2,3.3){\footnotesize{$(a,$-$b$+$1)$}}
\put(21.2,4.3){\footnotesize{$(a,$-$b$+$2)$}}
\put(14.6,8.3){\footnotesize{$(a$-$b,0)$}}
\put(16.0,8.3){\footnotesize{$(a$-$b$+$1,0)$}}
\put(0.1,2.3){$\mathcal{B}_{1,b}^{(r)}$}
\put(3.3,3.8){$\mathcal{B}_{b+1,b}^{(r)}$}
\put(5.3,4.6){$\mathcal{B}_{a,b}^{(r)}$}
\put(6.2,5.1){$\mathcal{A}_{a,b}^{(r)}$}
\put(7.2,5.6){$\mathcal{A}_{a,b-1}^{(r)}$}
\put(12.5,2.6){$\mathcal{B}_{1,b}^{(r)}$}
\put(17.7,5.0){$\mathcal{B}_{a,b}^{(r)}$}
\put(18.7,5.6){$\mathcal{A}_{a,b}^{(r)}$}
\put(19.7,6.1){$\mathcal{A}_{a,b-1}^{(r)}$}
\end{picture}
\caption{Case $a>b$}
\end{figure}
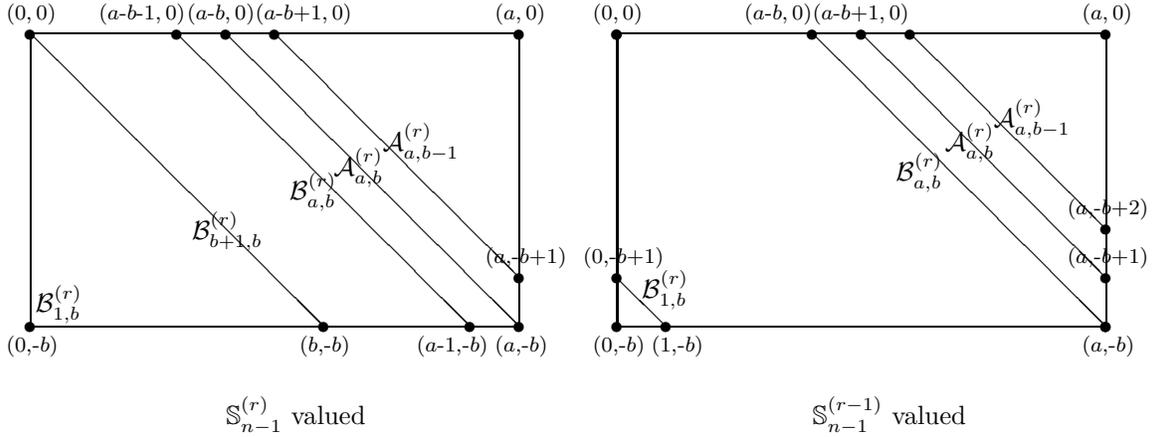

\unitlength 0.65cm
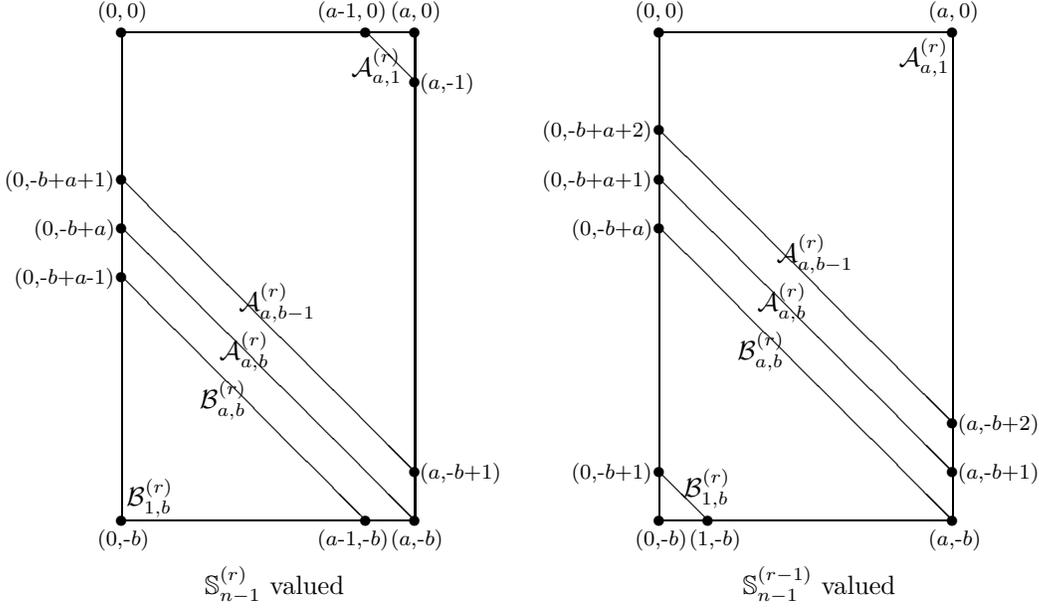
\begin{figure}[h]
\begin{picture}(17,12.5)(-3.1,0)
\put(-0.15,1.85){$\bullet$}
\put(0,2.0){\line(1,0){6}}
\put(4.85,1.85){$\bullet$}
\put(5.85,1.85){$\bullet$}
\put(6,2.0){\line(0,1){10}}
\put(5.85,2.85){$\bullet$}
\put(5.85,10.85){$\bullet$}
\put(5.85,11.85){$\bullet$}
\put(6,12.0){\line(-1,0){6}}
\put(4.85,11.85){$\bullet$}
\put(-0.15,11.85){$\bullet$}
\put(0,12.0){\line(0,-1){10}}
\put(-0.15,6.85){$\bullet$}
\put(-0.15,7.85){$\bullet$}
\put(-0.15,8.85){$\bullet$}
\put(0,7.0){\line(1,-1){5}}
\put(0,8.0){\line(1,-1){6}}
\put(0,9.0){\line(1,-1){6}}
\put(5,12.0){\line(1,-1){1}}

\put(10.85,1.85){$\bullet$}
\put(11,2.0){\line(1,0){6}}
\put(11.85,1.85){$\bullet$}
\put(16.85,1.85){$\bullet$}
\put(17,2.0){\line(0,1){10}}
\put(16.85,2.85){$\bullet$}
\put(16.85,3.85){$\bullet$}
\put(16.85,11.85){$\bullet$}
\put(17,12.0){\line(-1,0){6}}
\put(10.85,11.85){$\bullet$}
\put(11,12.0){\line(0,-1){10}}
\put(10.85,9.85){$\bullet$}
\put(10.85,8.85){$\bullet$}
\put(10.85,7.85){$\bullet$}
\put(10.85,2.85){$\bullet$}
\put(11,10.0){\line(1,-1){6}}
\put(11,9.0){\line(1,-1){6}}
\put(11,8.0){\line(1,-1){6}}
\put(11,3.0){\line(1,-1){1}}

\put(1.7,0.5){$\mS_{n-1}^{(r)}$ valued}
\put(12.7,0.5){$\mS_{n-1}^{(r-1)}$ valued}

\put(-0.5,12.3){\footnotesize{$(0,0)$}}
\put(-0.5,1.5){\footnotesize{$(0,$-$b)$}}
\put(5.5,12.3){\footnotesize{$(a,0)$}}
\put(4.0,12.3){\footnotesize{$(a$-$1,0)$}}
\put(5.5,1.5){\footnotesize{$(a,$-$b)$}}
\put(4.0,1.5){\footnotesize{$(a$-$1,$-$b)$}}
\put(10.5,12.3){\footnotesize{$(0,0)$}}
\put(10.5,1.5){\footnotesize{$(0,$-$b)$}}
\put(11.6,1.5){\footnotesize{$(1,$-$b)$}}
\put(16.5,12.3){\footnotesize{$(a,0)$}}
\put(16.5,1.5){\footnotesize{$(a,$-$b)$}}
\put(-2.2,6.85){\footnotesize{$(0,$-$b$+$a$-$1)$}}
\put(-1.8,7.85){\footnotesize{$(0,$-$b$+$a)$}}
\put(-2.4,8.85){\footnotesize{$(0,$-$b$+$a$+$1)$}}
\put(6.1,10.85){\footnotesize{$(a,$-$1)$}}
\put(6.1,2.85){\footnotesize{$(a,$-$b$+$1)$}}
\put(17.1,2.85){\footnotesize{$(a,$-$b$+$1)$}}
\put(17.1,3.85){\footnotesize{$(a,$-$b$+$2)$}}
\put(9.2,7.85){\footnotesize{$(0,$-$b$+$a)$}}
\put(8.6,8.85){\footnotesize{$(0,$-$b$+$a$+$1)$}}
\put(8.6,9.85){\footnotesize{$(0,$-$b$+$a$+$2)$}}
\put(9.2,2.85){\footnotesize{$(0,$-$b$+$1)$}}

\put(0.1,2.3){$\mathcal{B}_{1,b}^{(r)}$}
\put(1.6,4.3){$\mathcal{B}_{a,b}^{(r)}$}
\put(2.0,5.3){$\mathcal{A}_{a,b}^{(r)}$}
\put(2.4,6.3){$\mathcal{A}_{a,b-1}^{(r)}$}
\put(4.7,11.1){$\mathcal{A}_{a,1}^{(r)}$}
\put(11.5,2.5){$\mathcal{B}_{1,b}^{(r)}$}
\put(12.6,5.3){$\mathcal{B}_{a,b}^{(r)}$}
\put(13.0,6.3){$\mathcal{A}_{a,b}^{(r)}$}
\put(13.4,7.3){$\mathcal{A}_{a,b-1}^{(r)}$}
\put(15.9,11.3){$\mathcal{A}_{a,1}^{(r)}$}
\end{picture}
\caption{Case $a<b$}
\end{figure}


\section{Orthogonal bases of $\mcM_{a,b}^{(r)}$ in complex dimension $n=1$.}


In the complex plane ($n=1$) the self--adjoint idempotent $I$ is given by $I = \gf_1 \gfd_1$, and spinor space $\mS_1$ is the direct sum of two homogeneous parts $\mS_1^{(0)} = {\rm span}_{\mC}\{1\}I$ and $\mS_1^{(1)} = {\rm span}_{\mC}\{\gfd_1\}I$. For functions $F: \mC \longrightarrow \mS_1^{(0)}$ the equation $\upzd F = \gf_1 \p_{\olz_1} F = 0$ is trivially fulfilled, and $\upz F = \gfd_1 \p_{z_1} F = 0$ means that $F$ is anti-holomorphic in the variable $\olz_1$, forcing the first of the bidegrees of the homogeneous Hermitean monogenic polynomials to be zero. Each of the U$(1)$-modules $\mcM_{0,k}^{(0)}, k=0, 1, 2, \ldots$, is one-dimensional; it has highest weight $(-k)$ and it is spanned by the standard polynomial $P_{(-k)} = \frac{\olz_1^k}{k!} I$, which shows the differentiation property $\p_{\olz_1} P_{(-k)} = P_{(-k+1)}$ and hence also $(\p_{\olz_1})^k P_{(-k)} = P_{(0)} = I$.\\[-2mm]

On the other hand, for functions $G: \mC \longrightarrow \mS_1^{(1)}$, it is the equation $\upz G = \gfd_1 \p_{z_1} G = 0$ which is trivially fulfilled, while the equation $\upzd G = \gfd_1 \p_{\olz_1} G = 0$ means that $G$ is holomorphic in the variable $z_1$. For holomorphic homogeneous polynomials the second of the bidegrees must be zero. Each of the U$(1)$-modules $\mcM_{k,0}^{(1)}, k=0, 1, 2, \ldots$, is one-dimensional; it has highest weight $(k+1)$ and it is spanned by the standard polynomial $P_{(k+1)} = \frac{z_1^k}{k!} \gfd_1$, showing the differentiation property  $\p_{z_1} P_{(k+1)} = P_{(k)}$ and hence also $(\p_{z_1})^k P_{(k+1)} = P_{(1)} = \gfd_1$.


\section{Orthogonal bases of $\mcM_{0,b}^{(0)}(\mC^n)$.}\label{OGB0b}


We have already pointed out that the cases where $r=0$ and $r=n$ are to be treated separately from the general approach. In this section we consider the case $r=0$. \\[-2mm]

For functions $F$ taking their values in $\mS_n^{(0)} = {\rm span}_{\mC}\{ 1  \} I, I = \gf_1 \gfd_1 \cdots \gf_n\gfd_n$, the equation $\upzd F = \sum_{j=1}^n \gf_j \p_{\olz_j} F = 0$ is trivially satisfied. The equation $\upz F = \sum_{j=1}^n \gfd_j \p_{z_j} F = 0$ means that $F$ is anti--holomorphic in the variables $(\olz_1, \olz_2, \ldots, \olz_n)$, forcing the first of the bidegrees of the anti--holomorphic homogneous polynomials considered to be zero. This implies that the U$(n)$--module $\mcM_{0,b}^{(0)}(\mC^n)$ consists of all homogeneous polynomials of degree $b$ in the variables $(\olz_1, \olz_2, \ldots, \olz_n)$. The dimension of $\mcM_{0,b}^{(0)}(\mC^n)$ is $m_{0,b}^{(0)}(\mC^n) = \frac{(n+b-1)!}{(n-1)!b!}$ and its highest weight is $\underbrace{(0,0,\ldots,0,-b)}_{n}$. \\[-2mm]

The, very well--known, basis polynomials of $\mcM_{0,b}^{(0)}(\mC^n)$ are obtained through CK--extension from the initial polynomials in the spaces $\mcA_{0,b-j}^{(0)}, j=0,\ldots,b$, which are nothing else but the U$(n-1)$--modules $\mcM_{0,b-j}^{(0)}(\mC^{n-1})$ of dimension $m_{0,b-j}^{(0)}(\mC^{n-1}) = \frac{(n+b-j-2)!}{(n-2)!(b-j)!}$ and with highest weight $\underbrace{(0,0,\ldots,-b+j)}_{n-1} \prec \underbrace{(0,0,\ldots,0,-b)}_{n}$. Note that indeed 
$$
\sum_{j=0}^{b} m_{0,b-j}^{(0)}(\mC^{n-1}) = \sum_{j=0}^{b} \frac{(n+b-j-2)!}{(n-2)!(b-j)!} = \frac{(n+b-1)!}{(n-1)!b!} = m_{0,b}^{(0)}(\mC^n).
$$
Proceeding in the same way we will end up with the one--dimensional U$(1)$--modules $\mcM_{0,k}^{(0)}$, $k=0,\ldots,b$, with highest weights $(-k)$, as prescribed by the Gel'fand--Tsetlin procedure. A basis for $\mcM_{0,b}^{(0)}(\mC^n)$ is thus given, up to a multiplicative constant, by the polynomials
\begin{equation*}
P_{\underbrace{(0,\ldots,0,-b)}_{n}\underbrace{(0,\ldots,0,-b+j_{n-1})}_{n-1} \ldots \underbrace{(0,-b+j_{n-1}+ \cdots + j_2)}_{2} \underbrace{(-b+j_{n-1}+\cdots+j_1)}_{1}} =
\end{equation*}
$$
\frac{\olz_n^{j_{n-1}}}{j_{n-1}!} \frac{\olz_{n-1}^{j_{n-2}}}{j_{n-2}!} \cdots \frac{\olz_2^{j_1}}{j_1!} \frac{\olz_1^{(b-j_{n-1}- \cdots -j_1)}}{(b-j_{n-1}-\cdots-j_1)!} \ I
$$
where $j_{n-1}=0,\ldots,b$, $j_{n-2}=0,\ldots,b-j_{n-1},$ \ldots, $j_1=0,\ldots,b-j_{n-1}- \cdots -j_2$.\\[-2mm]

As an illustration let us consider the low dimensional cases $n=1, 2, 3$. For $n=1$, we have already seen in the foregoing section that the U$(1)$--module $\mcM_{0,b}^{(0)}(\mC)$ is spanned by the basis polynomial
$$
P_{(-b)} = \frac{\olz_1^b}{b!} I
$$
which shows the differentiation properties
\begin{eqnarray}\label{appell1b}
\p_{\olz_1} P_{(-b)} & = & P_{(-b+1)}\\
\p_{z_1} P_{(-b)} & = & 0
\end{eqnarray}

\vspace*{2mm}
For $n=2$, the U$(2)$--module $\mcM_{0,b}^{(0)}(\mC^2)$ has dimension $(b+1)$ and highest weight $(0,-b)$. Its basis stems, via CK-extension, from the initial data spaces $\mcA_{0,b-j_1}^{(0)}, j_1=0,\ldots,b$, which are nothing else but $\mcM_{0,b-j_1}^{(0)}(\mC) = {\rm span}_{\mC}\{ P_{(-b+j_1)} = \frac{\olz_1^{b-j_1}}{(b-j_1)!} I \}$, with corresponding highest weights $(-b+j_1), j_1=0,\ldots,b$. It follows that an orthogonal basis for $\mcM_{0,b}^{(0)}(\mC^2)$ is given by the polynomials
\begin{equation}\label{pol0b}
P_{(0,-b)(-b+j_1)} = \frac{\olz_2^{j_1}}{j_1!} \frac{\olz_1^{b-j_1}}{(b-j_1)!} I, \; j_1 = 0, \ldots, b
\end{equation}
showing the differentiation properties 
\begin{eqnarray}\label{appell0b}
\p_{\olz_2}P_{(0,-b)(-b+j_1)} & = & P_{(0,-b+1)(-b+j_1)}, \; \; \; \; \, j_1 = 1, \ldots, b\\
\p_{\olz_1}P_{(0,-b)(-b+j_1)} & = & P_{(0,-b+1)(-b+j_1+1)}, \; j_1 = 0, \ldots, b-1\nonumber
\end{eqnarray}

\vspace*{2mm}
For $n=3$, the U$(3)$--module $\mcM_{0,b}^{(0)}(\mC^3)$ has dimension $m_{0,b}^{(0)}(\mC^3) = \frac{(b+1)(b+2)}{2}$ and highest weight $(0,0,-b)$. Its basis is obtained through CK--extension of the initial polynomials contained in the spaces $\mcA_{0,b-j_2}^{(0)}$, $j_2 = 0, \ldots, b$. More explicitly:
\begin{eqnarray*}
\mcA_{0,0}^{(0)} \cong \mcM_{0,0}^{(0)}(\mC^2) & \stackrel{CK}{\longrightarrow} & P_{(0,0,-b)(0,0)(0)} = \frac{\olz_3^b}{b!} I \\
& \vdots & \\
\mcA_{0,b-j_2}^{(0)} \cong \mcM_{0,b-j_2}^{(0)}(\mC^2) & \stackrel{CK}{\longrightarrow} & P_{(0,0,-b)(0,-b+j_2)(-b+j_2+j_1)} = \frac{\olz_3^{j_2}}{j_2!} \frac{\olz_2^{j_1}}{j_1!} \frac{\olz_1^{(b-j_2-j_1)}}{(b-j_2-j_1)!}\ I, \\
& &  j_1 = 0, \ldots, b-j_2 \\
& \vdots & \\
\mcA_{0,b}^{(0)} \cong \mcM_{0,b}^{(0)}(\mC^2)& \stackrel{CK}{\longrightarrow} & P_{(0,0,-b)(0,-b)(-b+j_1)} =  \frac{\olz_2^{j_1}}{j_1!} \frac{\olz_1^{(b-j_1)}}{(b-j_1)!}\ I, \\
&&  j_1 = 0, \ldots, b
\end{eqnarray*}
Note the differentiation properties
\begin{eqnarray*}
\p_{\olz_3} P_{(0,0,-b)(0,-b+j_2)(-b+j_2+j_1)} &=& P_{(0,0,-b+1)(0,-b+j_2)(-b+j_2+j_1)}, \\
& &  j_2=1,\ldots,b, \; j_1=0,\ldots,b-j_2\\[2mm]
\p_{\olz_2} P_{(0,0,-b)(0,-b+j_2)(-b+j_2+j_1)} &=& P_{(0,0,-b+1)(0,-b+j_2+1)(-b+j_2+j_1)}, \\
& &  j_2=0,\ldots,b-1, \; j_1=1,\ldots,b-j_2\\[2mm]
\p_{\olz_1} P_{(0,0,-b)(0,-b+j_2)(-b+j_2+j_1)} &=& P_{(0,0,-b+1)(0,-b+j_2+1)(-b+j_2+j_1+1)}, \\
& & j_2=0,\ldots,b-1, \; j_1=0,\ldots,b-j_2-1
\end{eqnarray*}


\section{Orthogonal bases of $\mcM_{a,0}^{(n)}(\mC^n)$.}\label{OGBa0}


In this section we consider the exceptional case $r=n$. \\[-2mm]

For functions $F$ taking their values in $\mS_n^{(n)} = {\rm span}_{\mC}\{ \gfd_1 \gfd_2 \cdots \gfd_n  \}I = {\rm span}_{\mC}\{ \gfd_1 \gfd_2 \cdots \gfd_n  \}$, the equation $\upz F = \sum_{j=1}^n \gfd_j \p_{z_j} F = 0$ is trivially satisfied. The equation $\upzd F = \sum_{j=1}^n \gf_j \p_{\olz_j} F = 0$ means that $F$ is holomorphic in the variables $(z_1, z_2, \ldots, z_n)$, forcing the second of the bidegrees of the holomorphic homogeneous polynomials considered to be zero. This implies that the U$(n)$--module $\mcM_{a,0}^{(n)}(\mC^n)$ consists of all homogeneous polynomials of degree $a$ in the variables $(z_1, z_2, \ldots, z_n)$. The dimension of $\mcM_{a,0}^{(n)}(\mC^n)$ is $m_{a,0}^{(n)}(\mC^n) = \frac{(n+a-1)!}{(n-1)!a!}$ and its highest weight is $\underbrace{(a+1,1,\ldots,1)}_{n}$. \\[-2mm]

The, very well--known, basis polynomials of $\mcM_{a,0}^{(n)}(\mC^n)$ are obtained through CK--extension from the initial polynomials in the spaces $\mcB_{a-i,0}^{(n)}, i=0,\ldots,a$, which are nothing else but shifted versions (with embedding factor $\gfd_n$) of the U$(n-1)$--modules $\mcM_{a-i,0}^{(n)}(\mC^{n-1})$ of dimension $m_{a-i,0}^{(n)}(\mC^{n-1}) = \frac{(n+a-i-2)!}{(n-2)!(a-i)!}$ and with highest weight $\underbrace{(a+1-i,1,\ldots,1)}_{n-1} \prec \underbrace{(a+1,1,\ldots,1)}_{n}$. Note that indeed 
$$
\sum_{i=0}^{a} m_{a-i,0}^{(n)}(\mC^{n-1}) = \sum_{i=0}^{a} \frac{(n+a-i-2)!}{(n-2)!(a-i)!} = \frac{(n+a-1)!}{(n-1)!a!} = m_{a,0}^{(n)}(\mC^n).
$$
Proceeding in the same way we will end up with the one--dimensional U$(1)$--modules $\gfd_n \, \mcM_{k,0}^{(n)}, k=0,\ldots,a$, with highest weights $(k+1)$, as prescribed by the Gel'fand--Tsetlin procedure. A basis for $\mcM_{a,0}^{(n)}(\mC^n)$ is thus given by the polynomials
\begin{equation*}
P_{\underbrace{(a+1,1,\ldots,1)}_{n} \underbrace{(a+1-i_{n-1},1,\ldots,1)}_{n-1} \ldots \underbrace{(a+1-i_{n-1}-\cdots-i_2,1)}_{2} \underbrace{(a+1-i_{n-1}-\cdots-i_1)}_{1}} =
\end{equation*}
$$
\frac{z_n^{i_{n-1}}}{i_{n-1}!} \frac{z_{n-1}^{i_{n-2}}}{i_{n-2}!} \cdots \frac{z_2^{i_1}}{i_1!} \frac{z_1^{(a-i_{n-1}- \cdots -i_1)}}{(a-i_{n-1}-\cdots-i_1)!} \ \gfd_1 \cdots \gfd_n
$$
where $i_{n-1}=0,\ldots,a$, $i_{n-2}=0,\ldots,a-i_{n-1}$, \ldots , $i_1=0,\ldots,a-i_{n-1}- \cdots -i_2$.\\[-2mm]

As an illustration let us again consider the low dimensional cases $n=1, 2, 3$. For $n=1$, we have already seen in Section 4 that the U$(1)$--module $\mcM_{a,0}^{(1)}(\mC)$ is spanned by the basis polynomial
$$
P_{(a+1)} = \frac{z_1^a}{a!} \gfd_1
$$
showing the differentiation properties
\begin{eqnarray}\label{appell1a}
\p_{\olz_1} P_{(a+1)} & = & 0\\
\p_{z_1} P_{(a+1)} & = & P_{(a)}
\end{eqnarray}

\vspace{2mm}
For $n=2$, the U$(2)$--module $\mcM_{a,0}^{(2)}(\mC^2)$ has dimension $(a+1)$ and highest weight $(a+1,1)$. Its basis stems, via CK-extension, from the initial data spaces $\mcB_{a-i_1,0}^{(2)}, i_1=0,\ldots,a$, which are nothing else but $\gfd_2 \ \mcM_{a-i,0}^{(1)}(\mC)$, 
with corresponding highest weights $(a+1-i_1), i_1=0,\ldots,a$. 
As we have that
$$
\mcM_{a-i,0}^{(1)}(\mC) = {\rm span}_{\mC}\{ P_{(a+1-i_1)} = \frac{z_1^{a-i_1}}{(a-i_1)!} \gfd_1 I\}
$$
it follows that an orthogonal basis for $\mcM_{a,0}^{(2)}(\mC^2)$ is given by the polynomials
\begin{equation}\label{pola0}
P_{(a+1,1)(a+1-i_1)} = \frac{z_2^{i_1}}{i_1!} \frac{z_1^{a-i_1}}{(a-i_1)!} \gfd_1 \gfd_2 , \; i_1 = 0, \ldots, a
\end{equation}
showing the differentiation properties 
\begin{eqnarray}\label{appella0}
\p_{z_2}P_{(a+1,1)(a+1-i_1)} & = & P_{(a,1)(a+1-i_1)}, \;   i_1 = 1, \ldots, a\\
\p_{z_1}P_{(a+1,1)(a+1-i_1)} & = & P_{(a,1)(a-i_1)}, \; \; \; \; \, i_1 = 0, \ldots, a-1\nonumber
\end{eqnarray}

For $n=3$, the U$(3)$--module $\mcM_{a,0}^{(3)}(\mC^3)$ has dimension $m_{a,0}^{(3)}(\mC^3) = \frac{(a+1)(a+2)}{2}$ and highest weight $(a+1,1,1)$. Its orthogonal basis is obtained through CK--extension of the initial data spaces $\mcB_{a-i_2,0}^{(3)}, i_2 = 0, \ldots, a$, more explicitly:
\begin{eqnarray*}
\mcB_{a,0}^{(3)} \cong \gfd_3 \ \mcM_{a,0}^{(2)}(\mC^2) & \stackrel{CK}{\longrightarrow} & P_{(a+1,1,1)(a+1,1)(a+1-i_1)} = \frac{z_2^{i_1}}{i_1!}  \frac{z_1^{a-i_1}}{(a-i_1)!} \ \gfd_1 \gfd_2 \gfd_3 , \\
& & i_1=0,\ldots,a \\
& \vdots & \\
\mcB_{a-i_2,0}^{(3)} \cong \gfd_3 \ \mcM_{a-i_2,0}^{(2)}(\mC^2) & \stackrel{CK}{\longrightarrow} & P_{(a+1,1,1)(a+1-i_2,1)(a+1-i_2-i_1)} = \frac{z_3^{i_2}}{i_2!} \frac{z_2^{i_1}}{i_1!} \frac{z_1^{(a-i_2-i_1)}}{(a-i_2-i_1)!} \ \gfd_1 \gfd_2 \gfd_3 , \\
& & i_1 = 0, \ldots, a-i_2 \\
 & \vdots & \\
\mcB_{0,0}^{(3)} \cong \gfd_3 \ \mcM_{0,0}^{(2)}(\mC^2) & \stackrel{CK}{\longrightarrow} & P_{(a+1,1,1)(1,1)(1)} =  \frac{z_3^{a}}{a!} \ \gfd_1 \gfd_2 \gfd_3 
\end{eqnarray*}
These basis polynomials show the differentiation properties
\begin{eqnarray*}
\p_{z_3} P_{(a+1,1,1)(a+1-i_2,1)(a+1-i_2-i_1)} &=& P_{(a,1,1)(a+1-i_2,1)(a+1-i_2-i_1)}, \\
& &  i_2=1,\ldots,a, \; i_1=0,\ldots,a-i_2 \\[2mm]
\p_{z_2} P_{(a+1,1,1)(a+1-i_2,1)(a+1-i_2-i_1)} &=& P_{(a,1,1)(a-i_2,1)(a+1-i_2-i_1)}, \\
& &  i_2=0,\ldots,a-1, \;  i_1=1,\ldots,a-i_2 \\[2mm]
\p_{z_1} P_{(a+1,1,1)(a+1-i_2,1)(a+1-i_2-i_1)} &=& P_{(a,1,1)(a-i_2,1)(a-i_2-i_1)}, \\
& & i_2=0,\ldots,a-1,  \; i_1=0,\ldots,a-i_2-1
\end{eqnarray*}


\section{Orthogonal bases of $\mcM_{a,b}^{(r)}$ in complex dimension $n=2$}


In this section we will construct explicitly orthogonal bases for the spaces $\mcM_{a,b}^{(r)}(\mC^2)$ of spherical Hermitean monogenics $p(z_1, z_2, \olz_1, \olz_2)$ with values in the homogeneous parts $\mS_2^{(r)}$ of spinor space $\mS_2$. In this way we do not only provide an illustration of the algorithm devised in Section 3, but the bases themselves are important for applications in real dimension $4$. Spinor space $\mS_2$ decomposes into three homogeneous parts $\mS_2^{(0)} = {\rm span}_{\mC}\{1\}I$,  $\mS_2^{(1)} = {\rm span}_{\mC}\{\gfd_1, \gfd_2 \}I$, and   $\mS_2^{(2)} = {\rm span}_{\mC}\{\gfd_1 \gfd_2  \}I$. The cases $r=0$ and $r=2$ having been treated in the foregoing sections  for arbitrary complex dimension $n$, we concentrate on the case $r=1$.\\

The U($2$)-module $\mcM^{(1)}_{a,b}(\mC^2)$ has dimension $m_{a,b}^{(1)}(\mC^2) = a+b+2$, and highest weight $(a+1,-b)$. There are $(a+b+2)$ spaces of initial polynomials $\mcA^{(1)}_{a,b-j}$, $j=0,\ldots,b$ and $\mcB^{(1)}_{a-i,b}$, $i=0,\ldots,a$, each of them being one-dimensional, since $x_{a,b,1}=1$ and $y_{a,b,1}=1$. The general theory of the CK-extension procedure, see Theorem 3, predicts that the compatibility conditions imposed on these initial polynomials will be trivially fulfilled, which means that these initial polynomials are simply all homogeneous polynomials in the variables $z_1$ and $\olz_1$ of the appropriate bidegree. This is moreover confirmed by the Fischer decomposition. Indeed, taking into account that 
\begin{eqnarray*}
\tmcM_{0,b}^{(0)}(\mC) &=& {\rm span}_{\mC}\{P_{(-b)} = \frac{\olz_1^b}{b!} I \} \\
\tmcM_{a,0}^{(1)}(\mC) &=& {\rm span}_{\mC}\{P_{(a+1)} = \frac{z_1^a}{a!} \gfd_1 I \}
\end{eqnarray*}
by Theorem 4 we obtain
\begin{eqnarray*}
\mcA^{(1)}_{a,b-j} & = & \mbox{span}_{\mC} \left \{ \frac{{z_1}^a}{a!} \, \frac{{\olz_1}^{b-j}}{(b-j)!} \, \gfd_1 \, I \right \}, \qquad j=0,\ldots,b \\
\mcB^{(1)}_{a-i,b} & = & \mbox{span}_{\mC} \left \{ \frac{{z_1}^{a-i}}{(a-i)!} \, \frac{{\olz_1}^{b}}{b!} \, \gfd_2 \, I \right \}, \qquad i=0,\ldots,a
\end{eqnarray*}
Via CK-extension, each of the spaces of initial polynomials thus gives rise to exactly one Hermitean monogenic basis polynomial, establishing in this way an isomorphism between $\mcM^{(1)}_{a,b}(\mC^2)$, which is of dimension $(a+b+2)$, and the direct sum of one-dimensional subspaces
$$
\bigoplus_{j=0}^b \ \mcA_{a,b-j}^{(1)} \  \oplus \ \bigoplus_{i=0}^{a} \ \mcB_{a-i,b}^{(1)}
$$

Now let us check the branching rules for the space $\mcM^{(1)}_{a,b}(\mC^2)$ which we recall to have highest weight $\la = (a+1,-b)$. If $a > b$ then $\mcA_{a,b-j}^{(1)}$ is a shifted version of the U($1$)-module $\widetilde{\mcM}^{(1)}_{a-b+j,0}(\mC)$ with highest weight $(a-b+j+1)$, for all $j=0,\ldots,b$. At the same time $\mcB_{a-i,b}^{(1)}$ is either a shifted version of the U($1$)-module $\widetilde{\mcM}^{(1)}_{a-i-b-1,0}$ with highest weight $(a-i-b)$, for all $i=0,\ldots,a-b-1$, or a shifted version of the U$(1)$-module $\widetilde{\mcM}^{(0)}_{0,b-a+i}$ with highest weight $(a-i-b)$, for all $i=a-b,\ldots,a$. If $a < b$ then $\mcA_{a,b-j}^{(1)}$ is either a shifted version of the U($1$)-module $\widetilde{\mcM}^{(1)}_{a-b+j,0}(\mC)$ with highest weight $(a-b+j+1)$, for all $j=b-a,\ldots,b$, or it is a shifted version of the U$(1)$-module $\tmcM_{0,b-j-a-1}^{(0)}$, again with highest weight $(a+1+j-b)$, for all $j=0,\ldots,b-a-1$. At the same time $\mcB_{a-i,b}^{(1)}$ is a shifted version of the U($1$)-module $\widetilde{\mcM}^{(0)}_{0,b-a-i}$ with highest weight $(a-i-b)$, for all $i=0,\ldots,a$.\\[-2mm]

In both cases exactly all the highest weights $\mu \prec \la = (a+1,-b)$ are recovered, confirming that the CK-isomorphism yields the splitting of the irreducible U$(2)$-module $\mcM_{a,b}^{(1)}(\mC^2)$ into a direct sum of $(a+b+2)$ one-dimensional irreducible U$(1)$-modules $V_{\mu}, \mu = -b, \ldots, a+1$.\\[-2mm]

Finally an orthogonal basis for $\mcM^{(1)}_{a,b}(\mC^2)$ is explicitly constructed by applying the CK-extension to the initial polynomials contained in the spaces $\mcA_{a,b-j}^{(1)}, j=0,\ldots,b$, and $\mcB_{a-i,b}^{(1)}, i=0,\ldots,a$. The following closed form expressions are obtained for these basis polynomials:\\
(i) for $j=0,\ldots,b$ one has
\begin{equation}\label{pol1}
P_{(a+1,-b),(a+1-b+j)} = \sum_{k=0}^{\mbox{min}(a,b-j)} (-1)^{b-j-k} \, \frac{{z_2}^k}{k!} \, \frac{{\overline{z}_2}^{k+j}}{(k+j)!} \, \frac{{z_1}^{a-k}}{(a-k)!} \, \frac{{\overline{z}_1}^{b-j-k}}{(b-j-k)!} \, \gfd_1 \, I
\end{equation}
$$
 + \sum_{k=0}^{\mbox{min}(a,b-j-1)} (-1)^{b-j-k-1} \, \frac{{z_2}^k}{k!} \, \frac{{\overline{z}_2}^{k+j+1}}{(k+j+1)!} \, \frac{{z_1}^{a-k}}{(a-k)!} \, \frac{{\overline{z}_1}^{b-j-k-1}}{(b-j-k-1)!} \, \gfd_2 \, I
$$
(ii) for $i=0,\ldots,a$ one has
\begin{equation}\label{pol2}
P_{(a+1,-b),(a-b-i)} = \sum_{k=0}^{\mbox{min}(a-i,b)} (-1)^{b-k} \, \frac{{z_2}^{k+i}}{(k+i)!} \, \frac{{\overline{z}_2}^{k}}{k!} \, \frac{{z_1}^{a-i-k}}{(a-i-k)!} \, \frac{{\overline{z}_1}^{b-k}}{(b-k)!} \, \gfd_2 \, I
\end{equation}
$$
+ \sum_{k=0}^{\mbox{min}(a-i-1,b)} (-1)^{b-k} \, \frac{{z_2}^{k+i+1}}{(k+i+1)!} \, \frac{{\overline{z}_2}^{k}}{k!} \, \frac{{z_1}^{a-i-k-1}}{(a-i-k-1)!} \, \frac{{\overline{z}_1}^{b-k}}{(b-k)!} \, \gfd_1 \, I
$$
In the next section we study the properties of these basis polynomials under derivation with respect to the four variables $\olz_2, z_2, \olz_1, z_1$.


\section{The Appell property in complex dimension $n=2$}


In this section we show that in complex dimension $n=2$, the system of the constructed orthogonal bases of Hermitean monogenic polynomials possesses the Appell property with respect to all the variables, that is, by differentiating any basis element with respect to one of the variables $\olz_2$, $z_2$, $\olz_1$ or $z_1$, always a multiple of another basis element is obtained. This property is obvious for the $\mS_2^{(0)}$-- and $\mS_2^{(2)}$--valued basis polynomials, see \eqref{appell0b} in Section \ref{OGB0b} and \eqref{appella0} in Section \ref{OGBa0} respectively. For $\mS_2^{(1)}$--valued polynomials, the following result is obtained (see also \cite{rhofb}).

\begin{theorem} Let the $\mS_2^{(1)}$--valued basis polynomials $P_{(a+1,-b),(\mu)}$ be defined as in \eqref{pol1} and \eqref{pol2} above.
Then one has that
\begin{itemize}
\item[(i)] $\p_{\olz_2} P_{(a+1,-b),(\mu)}  =  P_{(a+1,-b+1),(\mu)}$;
\item[(ii)] $\p_{z_2} P_{(a+1,-b),(\mu)}  =  P_{(a,-b),(\mu)}$;
\item[(iii)] $\p_{\olz_1} P_{(a+1,-b),(\mu)}  =  - P_{(a+1,-b+1),(\mu+1)}$;
\item[(iv)] $\p_{z_1} P_{(a+1,-b),(\mu)}  = P_{(a,-b),(\mu-1)}$.
\end{itemize}
Here we have put $P_{(a+1,-b),(\mu)}=0$ unless $-b\leq \mu\leq a+1$.
\end{theorem}

\begin{remark}
The meaning of the Appell property formulated in Theorem 5 is the following. We can consider a finite dimensional subspace of Hermitean monogenic polynomials with bidegree of homogeneity bounded by fix constants $a$, $b$. On this subspace, each of the four derivatives is represented
in the constructed bases by a very simple nilpotent matrix. It has the property being a block matrix (with respect
to the decomposition given by the irreducible pieces of the considered subspace) where almost all blocks are zero matrices and where each non--zero block has the property that every column contains at most one nontrivial entry! In this way this Appell property can make numerical calculations
very effective.
\end{remark}

\begin{remark}
Note that the conditions (i) and (ii) of Theorem 5 are in fact properties which hold in any dimension $n$, since derivation w.r.t. the ''last variables''  $\overline{z}_n$ and $z_n$ is obviously $U(n-1)$-invariant. Thus we have that the derivative of any basis polynomial with respect to the variable $\overline{z}_n$ or $z_n$ equals a multiple of another basis element. In other words, in any dimension higher than $2$ the system of the considered orthogonal bases has still the Appell property with respect to the last variables, but not with respect to all variables. 
\end{remark}

 \section{Orthogonal bases of spinor valued spherical monogenics}

In \cite{fischer, damdeef}, spinor valued spherical monogenics in $\mR^{2n}$ are expressed in terms of Hermitean monogenic polynomials as follows.

\begin{theorem} 
Let $\mcM_{k}(\mR^{2n}, \mS_n)$ stand for the space of $\mS_n$-valued monogenic polynomials in $\mR^{2n}$ which are $k$-homogeneous.
Then, under the action of the group $U(n)$, the space $\mcM_{k}(\mR^{2n}, \mS_n)$ has a~multiplicity free irreducible decomposition
$$
\mcM_{k}(\mR^{2n}, \mS_n)=\left(\bigoplus_{a=0}^k\bigoplus_{r=0}^n \mcM_{a,k-a}^{(r)}(\mC^n)\right)\oplus
\left(\bigoplus_{a=0}^{k-1}\bigoplus_{r=1}^{n-1} 
\left( \frac{\uz}{a + n - r} + \frac{\uzd}{k - 1 - a + r}    \right)
\mcM_{a,k-1-a}^{(r)}(\mC^n)\right)
$$
In particular, this decomposition is orthogonal with respect to the Fischer (or any invariant) inner product.
\end{theorem}

As a~direct consequence of this decomposition, we can produce orthogonal bases of the spaces $\mcM_{k}(\mR^{2n}, \mS_n)$ composed of the constructed bases of homogeneous Hermitean monogenic polynomials. Indeed, we have the following result.

\begin{corollary} 
Let $\mcO_{a,b}^{(r)}$ denote a~GT basis of the space $\mcM_{a,b}^{(r)}(\mC^n).$ Then the set
$$
\mcO=\left(\bigcup_{a=0}^k\bigcup_{r=0}^n \mcO_{a,k-a}^{(r)}\right)\cup
\left(\bigcup_{a=0}^{k-1}\bigcup_{r=1}^{n-1} 
\left( \frac{\uz}{a + n - r} + \frac{\uzd}{k - 1 - a + r}    \right)
\mcO_{a,k-1-a}^{(r)}\right)
$$
is a~basis of the space $\mcM_{k}(\mR^{2n}, \mS_n)$ of spinor valued spherical monogenics which is orthogonal with respect to the Fischer (or any invariant) inner product.
\end{corollary}

\begin{remark}
For the monogenic polynomials of standard Clifford analysis, there exist already various constructions of orthogonal bases
(see, e.g., \cite{som, dss, cacao, cagubo, bock, lav, lasova, pvl, DLS_Hodge}). As pointed out in Remark 1, derivation with respect to each of the variables is represented, in a given basis, by a nilpotent matrix. The construction of the basis for spaces of standard spherical monogenics from Hermitean monogenic ones is clearly reducing the size of the block matrix due to the fact that the space of monogenics of a given degree splits into
a sum of many smaller $U(n)$-irreducible parts. Consequently, the size of the individual blocks is substantially reduced by such decomposition. On the other hand, it is clear that only a part of the above mentioned block will have the Appell property, while other blocks (with nontrivial embedding factors) will be more complicated and should be computed explicitly to see their form.
 \end{remark}


\section{Conclusion}


In this paper we have studied the construction of orthogonal bases for spaces of homogeneous Hermitean monogenic polynomials. We have taken advantage of the Gel'fand--Tsetlin approach, the fundamental symmetry group in Hermitean Clifford analysis being the unitary group. Two earlier obtained powerful theorems, viz. the Cauchy-Kovalevskaia extension theorem and the Fischer decomposition theorem were shown to be crucial in devising the corresponding algorithm which can be used in any dimension. For complex dimension $n=2$ the orthogonal bases were explicitly constructed and the Appell property with respect to all the variables for the system of the obtained orthogonal bases has been established. Note that for the explicit construction in higher dimension, we have developed a Maple programme, which will be discussed in the forthcoming paper \cite{bcdls}.


\section*{Acknowledgements}


R. L\'{a}vi\v{c}ka and V. Sou\v{c}ek acknowledge support by the institutional grant MSM 0021620839 and by grant GA CR 201/08/0397.



\end{document}